\documentclass[5p]{elsarticle}

\usepackage[hidelinks=true]{hyperref}

\journal{journal}

\usepackage{graphics} 
\usepackage{epsfig} 
\usepackage{amsmath}
\usepackage{amssymb}  
\usepackage{bbm}
\usepackage{url}
\usepackage[dvipsnames]{xcolor}
\usepackage{versions}

\includeversion{or-letters}
\excludeversion{techreport}

\DeclareMathOperator*{\argmin}{arg\,min}

\newcommand{\ds}{\displaystyle}
\newcommand{\defn}{\stackrel{\triangle}{=}}
\newcommand{\tth}{{}^{\text{th}}}
\newcommand{\dto}{\downarrow}

\newcommand{\smax}{S_{\max}}

\newcommand{\bE}{\mathbb{E}}

\newcommand{\bR}{\mathbb{R}}
\newcommand{\bZ}{\mathbb{Z}}

\newcommand{\cA}{\mathcal{A}}
\newcommand{\cB}{\mathcal{B}}
\newcommand{\cC}{\mathcal{C}}
\newcommand{\cM}{\mathcal{M}}
\newcommand{\cP}{\mathcal{P}}
\newcommand{\cR}{\mathcal{R}}
\newcommand{\cS}{\mathcal{S}}

\newcommand{\pd}{{(d)}}
\newcommand{\peps}{{(\epsilon)}}

\newcommand{\vq}{\boldsymbol{q}}
\newcommand{\va}{\boldsymbol{a}}
\newcommand{\vs}{\boldsymbol{s}}
\newcommand{\vu}{\boldsymbol{u}}

\newcommand{\vqbar}{\overline{\vq}}

\newcommand{\qbar}{\overline{q}}

\newcommand{\vx}{\boldsymbol{x}}
\newcommand{\vy}{\boldsymbol{y}}

\newcommand{\vzero}{\boldsymbol{0}}
\newcommand{\vone}{\boldsymbol{1}}
\newcommand{\vmu}{\boldsymbol{\mu}}
\newcommand{\vnu}{\boldsymbol{\nu}}

\newcommand{\E}[1]{\bE\left[#1\right]}
\newcommand{\Eq}[1]{\bE_{\vq}\left[#1\right]}
\newcommand{\Var}[1]{\text{Var}\left[#1\right]}
\newcommand{\ind}[1]{\mathbbm{1}_{\left\{#1 \right\}}}
\newcommand{\interior}[1]{\text{Int}\left(#1\right)}

\newcommand{\blue}{\color{blue}}
\renewcommand{\blue}{}

\newtheorem{thm}{Theorem}
\newtheorem{lem}[thm]{Lemma}
\newtheorem{claim}[thm]{Claim}
\newtheorem{proposition}[thm]{Proposition}
\newdefinition{definition}{Definition}
\newdefinition{rmk}{Remark}
\newproof{pf}{Proof}

\allowdisplaybreaks

\begin{document}

\begin{frontmatter}

\title{Throughput and Delay Optimality of Power-of-$d$ Choices in Inhomogeneous Load Balancing Systems}

\author{Daniela Hurtado-Lange${}^*$, Siva Theja Maguluri \\
d.hurtado@gatech.edu, siva.theja@gatech.edu}
\address{Department of Industrial and Systems Engineering, Georgia Institute of Technology, Atlanta, Georgia, United States\\
755 Ferst Drive NW, Atlanta, GA 30332}

%
%

\begin{abstract}
	It is well-known that the power-of-$d$ choices routing algorithm maximizes throughput and is heavy-traffic optimal in load balancing systems with homogeneous servers. However, if the servers are heterogeneous, throughput optimality does not hold in general. We find necessary and sufficient conditions for throughput optimality of power-of-$d$ choices when the servers are heterogeneous, and we prove that almost the same conditions are sufficient to show heavy-traffic optimality. {\blue Additionally, we generalize the sufficient condition for throughput optimality to a larger class of routing policies.}
\end{abstract}


\end{frontmatter}


\section{Introduction}\label{sec:introduction}

Load balancing systems are multi-server Stochastic Processing Networks (SPNs) in which there is a single stream of job arrivals. A single dispatcher routes arrivals to one of the queues immediately after they enter the system and, after being routed, the jobs wait in the corresponding line until the assigned server can process them. The policy used by the dispatcher to route the jobs is called a routing algorithm, and an essential goal when designing routing algorithms is to balance the workload of the servers in a way that delay is minimized, and the stability region of the SPN is maximal. When a routing algorithm achieves maximal stability region, it is said to be throughput optimal. For a formal definition of throughput optimality in the case of a load balancing system, we refer the reader to Definition \ref{def:throughput}.

The most basic algorithm is random routing, under which new arrivals are routed to a queue selected uniformly at random. Advantages of this routing algorithm are that the dispatcher does not require any information about the servers' speed or queue length.  
However, it has been proved that it is not delay optimal and, if the servers are heterogeneous, the stability region of the load balancing system under random routing is not maximal \cite{Mukh15_etal_JSQ2_stability}.

A popular routing algorithm is Join the Shortest Queue (JSQ), under which the new arrivals are routed to the server with the least number of jobs in line. It has been proved that {\blue JSQ is optimal among policies that do not know job durations, under several optimality criteria}. 
For example, \cite{winston_JSQ_1977,weber1978optimal} proved that 
JSQ maximizes the number of customers that complete service by a given time $t$. In \cite{winston_JSQ_1977}, Poisson arrivals and exponential job sizes are assumed, whereas \cite{weber1978optimal} relaxes these assumptions. In \cite{ephremides1980simple}, it is shown that JSQ minimizes the total time needed to finish processing all the jobs that arrive by a fixed time $t$. All these consider a continuous time model in a general setting, i.e., without taking any asymptotic regime. In \cite{JSQ_HT_optimality} it is proved that JSQ minimizes delay in the heavy-traffic regime, i.e., when the arrival rate approaches the maximum capacity of the system. This characteristic of a policy is known as heavy-traffic optimality. More recently, \cite{atilla} showed that JSQ is both throughput and heavy-traffic optimal in the context of a load balancing system operating in discrete time. In this case, instead of proving that delay is minimized, the authors prove that the total number of jobs in the system is minimized.
Even though JSQ is optimal under multiple criteria, a drawback is that it requires the dispatcher to know all the queue lengths at any point of time.

Comparing JSQ to random routing suggests a trade-off between the expected delay and the amount of information required by the routing algorithms. A policy that can be considered to be in between them is the power-of-$d$ choices algorithm, where $d$ is an integer between 1 and the total number of servers $n$. Under this algorithm, $d$ servers are sampled uniformly at random and the new arrivals are routed to the server with the shortest queue among these. If $d=1$, then power-of-$d$ is the same as random routing, and if $d=n$ it is the same as JSQ. In the case of load balancing systems with identical servers, it has been proved that even if $d=2$, 
power-of-$d$ choices is throughput and heavy-traffic optimal \cite{magsriyin_itc12_journal}. It has also been shown that power-of-$d$ choices yields substantial improvement in the tail probabilities of the queue lengths in mean-field regime (i.e., when the number of servers increases to infinity) \cite{mitzenmacher_po2,mitzenmacher_po2_2}. 
Also, for small values of $d$, the amount of information required by the dispatcher to route new arrivals is significantly smaller than under JSQ.

A disadvantage of power-of-$d$ choices is that throughput and delay optimality have been proved only when the servers are identical. If the service rates are different, there are known counterexamples for throughput optimality \cite{Mukh15_etal_JSQ2_stability}. In other words, if the servers are different, power-of-$d$ may reduce the stability region of the load balancing system. If the dispatcher knows the service rates, throughput and delay optimality of a modified version of power-of-$d$ choices have been proved in \cite{chen2012po2_journal,mukhopadhyay2016-PO2-heterogeneous}. In this adaptation, the probability of sampling each server is proportional to its mean service rate. 
However, we are interested in studying the cases when service rates are unknown to the dispatcher. 

The primary contribution of this paper is the computation of necessary and sufficient conditions for throughput optimality of power-of-$d$ choices, that only depend on the mean service rate vector. 
Specifically, we characterize a polytope within which the service rate vectors should lie. In particular, if the servers are identical our conditions are satisfied. Our result formalizes the idea that, in order to have throughput optimality, all the queues need to be sampled frequently enough. Then, given that power-of-$d$ selects $d$ queues uniformly at random, our result implies that the service rates of different servers should be close to each other; but not necessarily equal. 

{\blue In \cite{FossCher_1998stability} the authors address a similar question. They study stability of a general load balancing system, and they obtain sufficient conditions for throughput optimality. However, they approach the problem from a different perspective, and they provide conditions that depend on the queue length processes. In this paper, we provide conditions that only depend on the service rates and the sampling scheme. Hence, our conditions are easier to check.}

The second contribution of this paper is the computation of the joint distribution of the scaled queue lengths in heavy-traffic. We show that, if the heterogeneous service rates lie in the interior of the polytope proposed for throughput optimality, the load balancing system operating under power-of-$d$ choices has the same limiting distribution as a load balancing system operating under JSQ. Therefore, our results imply that power-of-$d$ choices is heavy-traffic optimal. 

Heavy-traffic means that we analyze the system when it is loaded to its maximum capacity. In the limit, many systems behave as if their dimension was smaller, phenomenon known as State Space Collapse (SSC). For the heterogeneous load balancing system operating under power-of-$d$ choices we prove that, in the limit, the $n$-dimensional queueing system behaves as a one-dimensional system, i.e., a single server queue. Then, we use this result to find the joint distribution of queue lengths. We develop our analysis in discrete time (i.e., in a time slotted fashion), so we use the notion of SSC developed in \cite{atilla}. Then, we find the joint distribution of the queue lengths using the Moment Generating Function (MGF) method introduced in \cite{Hurtado_transform_method}. Heavy-traffic analysis of the load balancing system operating under power-of-$d$ choices has been done in the past, but only under the assumption of identical and independent servers \cite{magsriyin_itc12_journal}. To the best of our knowledge, we are the first ones to obtain the heavy-traffic behavior of this queueing system with heterogeneous servers, and without modifying the probability of sampling each server.

{\blue The third contribution of this paper is a sufficient condition for throughput optimality under a larger class of routing policies. Specifically, we consider the following generalization of power-of-$d$ choices. In power-of-$d$ choices, only sets of size $d$ are sampled, and all of them are observed with the same probability. In the last part of this paper, we consider a routing policy that selects any subset of servers with certain probability, and routes the arrivals to the server with the shortest queue in the set. Then, we prove sufficient conditions on the sampling probabilities for throughput optimality.  
}

The organization of this paper is as follows. In Section \ref{sec:model} we formally introduce a model for the load balancing system and power-of-$d$ choices algorithm; in Section \ref{sec:throughput-optimality} we prove necessary and sufficient conditions for throughput optimality of power-of-$d$ choices; in Section \ref{sec:heavy-traffic} we perform heavy-traffic analysis; in Section \ref{sec:generalization} we present the generalization; and in Section \ref{sec:proofs} we present details of the proofs of the previous sections.

\subsection{Notation}\label{subsec:notation}

Before establishing the details of our model we introduce our notation. We use $\bR$ and $\bZ$ to denote the set of real and integer numbers, respectively. We add a subscript $+$ to indicate nonnegativity, and a number in the superscript to denote vector spaces. 
For any number $n\in\bZ_+$, we use $[n]\defn \left\{i\in\bZ_+:1\leq i\leq n \right\}$ and for $d\in\bZ_+$ with $n\geq d$ we use $\binom{n}{d}$ to denote the binomial coefficient. We use bold letters to denote vectors and the same letter but not bold and with a subscript $i$ to denote its $i\tth$ element. 
Given a vector $\vx\in\bR^n$, the notation $x_{(i)}$ refers to the $i\tth$ smallest element of $\vx$. 
Given two vectors $\vx,\vy\in\bR^n$, we use $\langle\vx,\vy\rangle$ to denote dot product and $\|\vx\|$ to the Euclidean norm. Then, $\|\vx\|=\sqrt{\langle\vx,\vx\rangle}$. Given a set $\cC\subset \bR^n$, we use $\interior{\cC}$ to denote its interior.

If $X$ is a random variable, then $\E{X}$ is its expected value and $\Var{X}$ its variance. For an event $A$, the notation $\ind{A}$ is the indicator function of $A$. Additionally, we use the notation $\Eq{\,\cdot\,}\defn \E{\,\cdot\, | \vq(k)=\vq}$ for the conditional expectation on the vector of queue lengths in time slot $k$. 

For any function $V:\bZ_+^n\to \bR_+$ let
\begin{align*}
	\Delta V(\vq)\defn\big[V\big(\vq(k+1)\big)-V\big(\vq(k)\big) \big]\ind{\vq(k)=\vq}.
\end{align*}
Thus, $\Delta V(\vq)$ is a random variable that measures the amount of change in the value of $V$ in one step, starting from $\vq$. We refer to $\Delta V(\vq)$ as the drift of $V(\vq)$.

\section{Model}\label{sec:model}

We model the load balancing system in discrete time, i.e., in a time slotted fashion, and we use $k\in\bZ_+$ to index time. Consider a system with $n$ servers, each of them with an infinite buffer. Let $\vq(k)$ be the vector of queue lengths at the beginning of time slot $k$, i.e., for each $i\in[n]$, $q_i(k)$ is the number of jobs in queue $i$ at the beginning of time slot $k$ including the job in service, if any. There is a single stream of arrivals to the system, and a dispatcher routes all arrivals of each time slot to one of the queues, according to some routing policy. We assume the routing time is negligible. Let $\{a(k):k\in \bZ_+\}$ be a sequence of i.i.d. random variables such that $a(k)$ is the total number of arrivals in time slot $k$. The vector $\va(k)$ represents the number of jobs that arrive to each of the queues in time slot $k$ after routing. Then, if the dispatcher routes the arrivals to queue $i^*$, we have $a_{i^*}(k)=a(k)$ and $a_i(k)=0$ for all $i\neq i^*$. Let $\vs(k)$ be the potential service vector in time slot $k$, i.e., for each $i\in[n]$, $s_i(k)$ is the number of jobs that can be processed in queue $i$ in time slot $k$ if there are enough jobs in line. Let $\{\vs(k):k\in \bZ_+\}$ be a sequence of i.i.d. random vectors, which is independent of the arrival and queue length processes. The difference between potential and actual service is called unused service, and we use $\vu(k)$ to denote the vector of unused service in time slot $k$. Observe that $\vu(k)$ is a function of $\vq(k)$, $\va(k)$ and $\vs(k)$.

We assume arrivals and routing occur before service in each time slot. Then, the following equation describes the dynamics of the queues. For each $i\in[n]$ and each $k\in \bZ_+$,
\begin{align}\label{eq:q-dynamics}
q_i(k+1)=q_i(k)+a_i(k)-s_i(k)+u_i(k).
\end{align}
From \eqref{eq:q-dynamics}, observe that $\left\{\vq(k):k\in \bZ_+ \right\}$ is a Discrete Time Markov chain (DTMC). Also, for every $i\in[n]$,
\begin{align}\label{eq:qu}
q_i(k+1)u_i(k)=0 
\end{align}
because the unused service in queue $i$ is nonzero only if the potential service to that queue is larger than the number of jobs available to be served (queue length and arrivals). Therefore, if the unused service is nonzero, the queue is empty at the beginning of the next time slot.

{\blue We assume that the arrival and the potential service to each queue have finite second moment.} Let $\lambda\defn\E{a(1)}$, $\vmu\defn\E{\vs(1)}$ and $\mu_\Sigma\defn \sum_{i=1}^n \mu_i$. Without loss of generality, we assume the vector $\vmu$ is ordered from minimum to maximum, i.e., $\mu_i=\mu_{(i)}$ for all $i\in[n]$. Let $\sigma_a^2\defn\Var{a(1)}$ be the variance of the arrival process and $\Sigma_s$ the covariance matrix of $\vs(1)$. It is well known that the capacity region of the load balancing system is
\begin{align}\label{eq:capacity-region}
\cC \defn \left\{\lambda\in\bR_+:\lambda\leq \mu_{\Sigma} \right\},
\end{align}
i.e., for each $\lambda\in\interior{\cC}$, there exists a routing algorithm such that $\{\vq(k):k\in \bZ_+\}$ is positive recurrent, and if $\lambda\notin \cC$, then $\{\vq(k):k\in \bZ_+\}$ is not positive recurrent for any 
routing algorithm. A proof of this is presented in \cite{atilla}. 

In this paper we work with 
power-of-$d$ choices, also known as JSQ($d$). We briefly describe it below.

\begin{definition}\label{def:po-d}
	Fix $d\in[n]$. In each time slot, the power-of-$d$ choices algorithm selects $d$ queues uniformly at random, and then routes the arrivals to the shortest of these. Ties are broken at random. Formally, if queues $i_1,\ldots,i_d$ are selected uniformly at random, then the arrivals in time slot $k$ are routed to the $i^*\tth$ queue, where \newline $i^*\in\argmin_{i\in\{i_1,\ldots,i_d\}}\left\{q_i(k) \right\}$.
\end{definition}

Observe that power-of-$d$ choices algorithm does not require any information about arrival or service rates. It just requires observing the number of jobs at $d$ of the queues.

\section{Throughput optimality of power-of-$d$ choices}\label{sec:throughput-optimality}

In this section we state and prove the main theorem of this paper. Before presenting the result we formally define throughput optimality.

\begin{definition}\label{def:throughput}
	A routing algorithm $\cA$ is throughput optimal if the queue length process $\left\{\vq(k):k\in \bZ_+ \right\}$ of the load balancing system operating under $\cA$ is positive recurrent for all $\lambda\in\interior{\cC}$, where $\cC$ is defined in \eqref{eq:capacity-region}.
\end{definition}

Now we present the main theorem of this paper. 
\begin{techreport}
	Then, $x_{(1)}=\min_{i\in[n]} x_i$ and $x_{(n)}=\max_{i\in[n]}x_i$, for example.
\end{techreport}

\begin{thm}\label{thm:throughput-optimality}
	For any $d\in[n-1]$, define
	\begin{align}\label{eq:Md}
	\cM^\pd\defn \left\{\vmu\in\bR^n_+: \dfrac{\sum_{i=1}^j \mu_{(i)}}{\mu_\Sigma}\geq \dfrac{\binom{j}{d}}{\binom{n}{d}}\; \forall d\leq j\leq n-1 \right\}.
	\end{align} 
	Then, the power-of-$d$ choices algorithm is throughput optimal for the load balancing system described in Section \ref{sec:model} if and only if $\vmu\in\cM^\pd$.
\end{thm}

\begin{rmk}
	Observe that we can equivalently define $\cM^\pd$ for all $d\in[n]$ as follows
	\begin{align*}
	\cM^\pd\defn \left\{\vmu\in\bR^n_+: \dfrac{\sum_{i=1}^j \mu_{(i)}}{\mu_\Sigma}\geq \dfrac{\binom{j}{d}}{\binom{n}{d}}\; \forall j\in [n] \right\},
	\end{align*}
	where we use the convention $\binom{j}{d}=0$ if $j<d$. Here we only added redundant constraints to $\cM^\pd$, so we use the definition \eqref{eq:Md} to avoid confusion.
\end{rmk}

\begin{rmk}
	An interpretation of Theorem \ref{thm:throughput-optimality} is the following. In order for power-of-$d$ choices algorithm to be throughput optimal, faster servers should be sampled sufficiently often. If this does not happen, it leads to the counter example in \cite{Mukh15_etal_JSQ2_stability}. Equation \eqref{eq:Md} characterizes the amount of imbalance between service rates that power-of-$d$ choices can tolerate. Note that, when the number of servers is fixed, as $d$ increases, power-of-$d$ choices can tolerate more imbalance because the right hand side of \eqref{eq:Md} becomes smaller. If $d=1$, which corresponds to random routing, the set $\cM^\pd$ is exactly the set of vectors where all the service rates are equal. In the other extreme case, when $d=n$, all the inequalities in \eqref{eq:Md} are redundant, and $\cM^\pd$ is the set of all nonnegative vectors. This fact is consistent with the throughput optimality of JSQ for any vector of service rates.
\end{rmk}

\begin{rmk}
	For $i\in[n]$, define $\nu_i \defn	\dfrac{\binom{i-1}{d-1}}{\binom{n}{d}}$, and let $\vnu$ be a vector with elements $\nu_i$.
	An equivalent characterization of $\cM^\pd$ is the set of all nonnegative vectors $\vmu$ such that $\frac{\vmu}{\mu_\Sigma}$ is majorized by $\vnu$. Majorization captures the notion of imbalance, and several equivalent characterizations can be found in \cite{majorization_book_MarshalOlkinArnold}. This notion has been used in the study of balls and bins models \cite{Azar_BallsAndBins}, and to prove optimality of routing and servicing algorithms \cite{Menich_OptimalityJSQ_majorization}. This notion also shows that for fixed $d$ and $n$, the vector $\vmu=\vnu$ is on the boundary of $\cM^\pd$.
\end{rmk}

{\blue \begin{rmk}
		Theorem \ref{thm:throughput-optimality} establishes that if $\vmu\notin \cM^\pd$, then the power-of-$d$ choices is not throughput optimal. In other words, if $\vmu\notin \cM^\pd$ there are some values of $\lambda\in\interior{\cC}$ for which $\left\{\vq(k):k\in\bZ_+ \right\}$ is not positive recurrent. In fact, if $\vmu\notin \cM^\pd$, the queue length process is positive recurrent only if $\lambda\in\interior{\overline{\cC}}$, where 
		\begin{align*}
			\overline{\cC} \defn \left\{\lambda\in\bR_+:\; \lambda\leq \dfrac{\binom{n}{d}}{\binom{j}{d}}\sum_{i=1}^j \mu_i \quad \forall d-1\leq j\leq n-1\right\}.
		\end{align*} 
		Observe that $\overline{\cC}\subsetneq\cC$ if $\vmu\notin \cM^\pd$, and $\cC=\overline{\cC}$ if $\vmu\in\cM^\pd$. We omit the proof of this remark, since it easily follows from the proof of Theorem \ref{thm:throughput-optimality}.
\end{rmk} }

In the proof of Theorem \ref{thm:throughput-optimality} we use Foster-Lyapunov theorem \cite[Theorem 3.3.7]{srikantleibook} and a certificate that a DTMC is not positive recurrent \cite[Theorem 3.3.10]{srikantleibook}. 
\begin{techreport}
	We state both of them in \ref{app:F-L.theorem} for completeness.
\end{techreport}

\begin{pf}[of Theorem \ref{thm:throughput-optimality}]
	Let $\epsilon\defn \mu_\Sigma-\lambda$, and observe that $\lambda\in\interior{\cC}$ if and only if $\epsilon\in(0,\mu_\Sigma)$. 
	We first prove that if $\vmu\in\cM^\pd$, then the power-of-$d$ choices algorithm is throughput optimal. To do that, we use Foster-Lyapunov theorem 
	\begin{techreport}
		(Theorem \ref{lemma:foster-lyapunov}) 
	\end{techreport}
	with Lyapunov function $V(\vq)=\left\|\vq \right\|^2$. We have
	\begin{techreport}
		\begin{align}
		& \Eq{\Delta V(\vq(k))} \nonumber \\
		=& \Eq{\left\|\vq(k+1) \right\|^2 - \left\|\vq(k) \right\|^2} \nonumber \\
		\stackrel{(a)}{=}& \Eq{\left\|\vq(k+1)-\vu(k) \right\|^2 + \left\|\vu(k)\right\|^2 \right. \nonumber \\
			&\left.+ 2\langle\vq(k+1)-\vu(k),\vu(k)\rangle - \left\|\vq(k) \right\|^2} \nonumber\\
		\stackrel{(b)}{=}& \Eq{\left\|\vq(k)+\va(k)-\vs(k) \right\|^2 - \left\|\vu(k) \right\|^2 - \left\| \vq(k)\right\|^2} \nonumber \\
		\stackrel{(c)}{\leq}& \Eq{\left\|\vq(k)+\va(k)-\vs(k) \right\|^2 - \left\| \vq(k)\right\|^2} \nonumber \\
		\stackrel{(d)}{=}& \Eq{\left\| \va(k)-\vs(k)\right\|^2} + 2\Eq{\langle\vq,\va(k)-\vs(k)\rangle}, \label{eq:foster-lyapunov-partial}
		\end{align}
		where $(a)$ holds after adding and subtracting $\vu(k)$ to the first term, and expanding the square; $(b)$ holds after using \eqref{eq:q-dynamics} and \eqref{eq:qu}, and reorganizing terms; $(c)$ holds because $\left\|\vu(k)\right\|^2\geq 0$; and $(d)$ holds after expanding the first square and reorganizing terms.
	\end{techreport}
	\begin{or-letters}
		\begin{align}
		& \Eq{\Delta V(\vq(k))} = \Eq{\left\|\vq(k+1) \right\|^2 - \left\|\vq(k) \right\|^2} \nonumber \\
		& \stackrel{(a)}{=} \Eq{\left\|\vq(k)+\va(k)-\vs(k) \right\|^2 - \left\|\vu(k) \right\|^2 - \left\| \vq(k)\right\|^2} \nonumber \\
		& \stackrel{(b)}{\leq}  \Eq{\left\| \va(k)-\vs(k)\right\|^2} + 2\Eq{\langle\vq,\va(k)-\vs(k)\rangle}, \label{eq:foster-lyapunov-partial}
		\end{align}
		where $(a)$ holds after few algebraic steps, using \eqref{eq:q-dynamics} and \eqref{eq:qu}; and $(b)$ holds because $\left\|\vu(k) \right\|^2\geq 0$ and after expanding the first term. 
	\end{or-letters}
	We analyze each of the terms in \eqref{eq:foster-lyapunov-partial} separately. 
	{\blue For the first term, we have
	\begin{align*}
	&\E{\left\|\va(k)-\vs(k) \right\|^2} \leq \E{\|\va(k)\|^2}+\E{\|\vs(k)\|^2} \\
	& \stackrel{(a)}{=} \E{a(k)^2} +\sum_{i=1}^n\E{s_i(k)^2} 
	\stackrel{(b)}{=} \lambda^2+\sigma_a^2+\sum_{i=1}^n\left(\mu_i^2+\sigma_{s_i}^2\right),
	\end{align*}
	where $(a)$ holds because all the arrivals in one time slot are routed to the same queue; and $(b)$ holds by definition of variance. Define $K_1\defn \lambda^2+\sigma_a^2+\sum_{i=1}^n\left(\mu_i^2+\sigma_{s_i}^2\right)$, and observe $K_1$ is a finite constant. Then, 
	\begin{align}\label{eq:foster-lyapunov-partial1}
	\Eq{\left\|\va(k)-\vs(k)\right\|^2}\leq K_1.
	\end{align}}
	
	{\blue Observe that the computation of the bound \eqref{eq:foster-lyapunov-partial1} does not use any properties of the routing algorithm. In other words, the bound \eqref{eq:foster-lyapunov-partial1} is valid for the load balancing system under any routing algorithm.} 
	
	To compute the second term of \eqref{eq:foster-lyapunov-partial}, we first compute $\Eq{\langle\vq,\va(k)\rangle}$. Recall that under power-of-$d$ choices, $d$ queues are chosen uniformly at random, and then the arrivals are sent to the shortest among them. Then, we have
	\begin{align}\label{eq:Ea}
	\Eq{\langle\vq,\va(k)\rangle}=& \lambda\sum_{i=1}^{n-d+1} q_{(i)}\dfrac{\binom{n-i}{d-1}}{\binom{n}{d}}
	\end{align}
	because there are $\binom{n-i}{d-1}$ ways to sample $d$ queues, and make sure that $q_{(i)}$ is the shortest; and 
	there are $\binom{n}{d}$ ways to sample $d$ queues uniformly at random.
	
	Let $\phi(i)$ be the index of the $i\tth$ shortest queue given $\vq(k)=\vq$. Then, since the potential service is independent of the queue lengths, the second term of \eqref{eq:foster-lyapunov-partial} is
	\begin{align}
	& \Eq{\langle\vq,\va(k)-\vs(k)\rangle} = \Eq{\langle\vq,\va(k)\rangle}- \langle\vq,\vmu\rangle \nonumber \\
	&= \sum_{i=1}^{n-d+1} q_{(i)}\left(\dfrac{\lambda\binom{n-i}{d-1}}{\binom{n}{d}} - \mu_{\phi(i)}\right)  - \sum_{i=n-d+2}^nq_{(i)} \mu_{\phi(i)}. \label{eq:alpha-partial}
	\end{align}
	Define
	\begin{align}\label{eq:def-alpha}
	\alpha_i\defn \begin{cases}
	\dfrac{\lambda\binom{n-i}{d-1}}{\binom{n}{d}} - \mu_{\phi(i)} & \text{, if }1\leq i\leq n-d+1 \\
	-\mu_{\phi(i)} & \text{, if }n-d+1< i\leq n.
	\end{cases} 
	\end{align}
	
	\begin{claim}\label{claim:alphas}
		The parameters $\alpha_i$ defined in \eqref{eq:def-alpha} satisfy 
		\begin{enumerate}
			\item\label{alpha:prop-n} $\alpha_n \leq -\mu_{(1)}$.
			
			\item\label{alpha:total-sum} $\ds\sum_{i=1}^n \alpha_i=-\epsilon$.
			
			\item\label{alpha:partial-sum} For any $j\in\bZ_+$ satisfying $2\leq j\leq n-1$,  we have $\ds \sum_{i=j}^n \alpha_i\leq -K_2$, where $K_2\defn \min\left\{\mu_{(1)},\frac{\epsilon}{\binom{n}{d}}\right\}$.
		\end{enumerate}
	\end{claim}
	
	We prove Claim \ref{claim:alphas} in Section \ref{app:proof-alphas}. Now we compute an upper bound for \eqref{eq:alpha-partial}. We obtain
	\begin{align}
	& \Eq{\langle\vq,\va(k)-\vs(k)\rangle}= \sum_{i=1}^n \alpha_i q_{(i)} \nonumber\\
	=& q_{(1)}\sum_{i=1}^n \alpha_i + \sum_{j=2}^n\left(\sum_{i=j}^n \alpha_i\right)\left(q_{(j)}-q_{(j-1)} \right) \nonumber \\
	\stackrel{(a)}{\leq}& -\epsilon q_{(1)} - K_2 \sum_{j=2}^n \left(q_{(j)}-q_{(j-1)}\right) \nonumber \\
	\stackrel{(b)}{=}& q_{(1)}\left(K_2-\epsilon\right) - K_2 q_{(n)} \stackrel{(c)}{\leq} - K_2 q_{(n)}, \label{eq:foster-lyapunov-partial2}
	\end{align}
	where $(a)$ holds by properties \ref{alpha:total-sum} and \ref{alpha:partial-sum} in Claim \ref{claim:alphas}; $(b)$ holds after solving the telescopic sum and rearranging terms; and $(c)$ holds because $K_2\leq \frac{\epsilon}{\binom{n}{d}}$ by definition, and $\binom{n}{d}\geq 1$. Using \eqref{eq:foster-lyapunov-partial1} and \eqref{eq:foster-lyapunov-partial2} in \eqref{eq:foster-lyapunov-partial} we obtain
	\begin{align*}
	\Eq{\Delta V(\vq(k))}\leq K_1 - 2K_2 q_{(n)}.
	\end{align*}
	
	\begin{or-letters}
		This inequality is sufficient to prove the conditions of Foster-Lyapunov theorem. 
	\end{or-letters}
	\begin{techreport}
		Defining
		\begin{align*}
		\cB\defn \left\{\vq\in\bR^n_+: \max_{i\in[n]}q_i \leq \dfrac{K_1+\xi}{2K_2} \right\},
		\end{align*}
		both of the conditions of Theorem \ref{lemma:foster-lyapunov} are satisfied.
	\end{techreport}
	Therefore, if $\vmu\in\cM^\pd$ then the power-of-$d$ choices algorithm is throughput optimal.
	
	Now we prove that if $\vmu\notin \cM^\pd$, then the power-of-$d$ choices algorithm is not throughput optimal. In other words, we prove that if $\vmu\notin \cM^\pd$, there exists $\lambda\in\interior{\cC}$ such that $\left\{\vq(k):k\in \bZ_+ \right\}$ is not positive recurrent. 
	
	First observe that if $\vmu\notin \cM^\pd$, there exists $j\in\bZ_+$ such that $d\leq j\leq n-1$ and $\frac{\sum_{i=1}^{j}\mu_i}{\mu_\Sigma}< \frac{\binom{j}{d}}{\binom{n}{d}}$. Let $j^*$ be smallest $j$ satisfying this condition, and $\delta_{j^*}>0$ satisfy
	\begin{align}\label{eq:def-delta-j}
	\dfrac{\sum_{i=1}^{j^*}\mu_i}{\mu_\Sigma} + \delta_{j^*}= \dfrac{\binom{j^*}{d}}{\binom{n}{d}}.
	\end{align}
	\begin{techreport}
		We use Lemma \ref{lemma:foster-lyapunov-not-pos-rec} with function $V_{j^*}(\vq)=\sum_{i=1}^{j^*}q_i$. We have
		\begin{align*}
		& \Eq{V_{j^*}(\vq(k+1))- V_{j^*}(\vq(k))} \\
		=& \sum_{i=1}^{j^*}\Eq{a_i(k) - s_i(k)+ u_i(k)} \\
		\stackrel{(a)}{\geq}& \sum_{i=1}^{j^*} \Eq{a_i(k)} - \sum_{i=1}^{j^*} \mu_i \\
		\stackrel{(b)}{\geq}& \sum_{i=1}^{j^*}\Eq{a_{\tilde{\phi}(i)}(k)} - \sum_{i=1}^{j^*} \mu_i \\
		\stackrel{(c)}{=}& \sum_{i=d}^{j^*}\lambda \dfrac{\binom{i-1}{d-1}}{\binom{n}{d}} - \mu_\Sigma\left(\dfrac{\binom{j^*}{d}}{\binom{n}{d}}-\delta_{j^*} \right) \\
		\stackrel{(d)}{=}& \mu_\Sigma \delta_{j^*} - \epsilon\dfrac{\binom{j^*}{d}}{\binom{n}{d}}
		\end{align*}
		where $(a)$ holds because $\E{s_i(k)}=\mu_i$ and $\E{u_i(k)}\geq 0$ for all $i\in[n]$; $(b)$ holds by letting $\tilde{\phi}(i)$ be the index of the $i\tth$ longest element of $\vq$, and because under power-of-$d$ choices the arrivals are routed to the shortest queue among the $d$ selected; $(c)$ holds by \eqref{eq:def-delta-j}, and because the arrivals are routed to the $i\tth$ longest queue only if the other $d-1$ selected queues are larger, and this happens with probability $\frac{\binom{i-1}{d-1}}{\binom{n}{d}}$ if $i\geq d$ and with probability 0 otherwise (similarly to the computation of \eqref{eq:Ea}); and $(d)$ holds because $\sum_{i=d}^{j^*}\binom{i-1}{d-1}=\binom{j^*}{d}$ and $\lambda=\mu_\Sigma -\epsilon$.
		
		This proves conditions (C1) and (C2) for $\epsilon>0$ satisfying 
		\begin{align*}
		\epsilon\leq \mu_\Sigma\min\left\{1, \delta_{j^*} \dfrac{\binom{j^*}{d}}{\binom{n}{d}} \right\}
		\end{align*}
		
		To prove condition (C3) observe
		\begin{align*}
		& \Eq{V_{j^*}(\vq(k+1))- V_{j^*}(\vq(k))} \\
		=& \sum_{i=1}^{j^*} \Eq{a_i(k)-\left(s_i(k)-u_i(k) \right)} \\
		\stackrel{(a)}{\leq}& \sum_{i=1}^{j^*} \Eq{a_i(k)} \\
		\stackrel{(b)}{\leq}& \sum_{i=1}^n \Eq{a_i(k)} \\
		\stackrel{(c)}{=}& \lambda
		\end{align*}
		where $(a)$ holds because $u_i(k)\leq s_i(k)$ with probability 1, by definition of unused service; $(b)$ holds because arrivals to each queue are a nonnegative random variable; and $(c)$ holds because $a(k)=\sum_{i=1}^n a_i(k)$ and $\lambda=\E{a(k)}$. Since $\lambda<\infty$, this proves condition $(C3)$.
	\end{techreport}
	\begin{or-letters}
		Using the Lyapunov function $V_{j^*}(\vq)=\sum_{i=1}^{j^*}q_i$, we obtain
		\begin{align*}
		& \Eq{\Delta V_{j^*}(\vq(k))} = \sum_{i=1}^{j^*}\Eq{a_i(k) - s_i(k)+ u_i(k)} \\
		& \stackrel{(a)}{\geq} \sum_{i=1}^{j^*} \Eq{a_i(k)} - \sum_{i=1}^{j^*} \mu_i \\
		& \stackrel{(b)}{\geq} \sum_{i=1}^{j^*}\Eq{a_{\tilde{\phi}(i)}(k)} - \sum_{i=1}^{j^*} \mu_i 
		 \stackrel{(c)}{=} \mu_\Sigma \delta_{j^*} - \epsilon\dfrac{\binom{j^*}{d}}{\binom{n}{d}},
		\end{align*}
		where $(a)$ holds because $\E{s_i(k)}=\mu_i$ and $\E{u_i(k)}\geq 0$ for all $i\in[n]$; $(b)$ holds by letting $\tilde{\phi}(i)$ be the index of the $i\tth$ longest element of $\vq$, and because under power-of-$d$ choices the arrivals are routed to the shortest queue among the $d$ selected; and $(c)$ holds computing $\Eq{a_{\tilde{\phi}(i)}(k)}$ similarly to \eqref{eq:Ea}, and reorganizing terms. 
		
		If $\epsilon>0$ satisfies $\epsilon\leq \mu_\Sigma\min\left\{1, \delta_{j^*} \frac{\binom{j^*}{d}}{\binom{n}{d}} \right\}$, then we have $\Eq{V_{j^*}(\vq(k+1))- V_{j^*}(\vq(k))}\geq 0$ for all $\vq\in\bR^n_+$. Additionally, we need to prove that $\Eq{\Delta V_{j^*}(\vq(k))}$ is bounded. We have
		\begin{align*}
		& \Eq{\Delta V_{j^*}(\vq(k))} 
		 \stackrel{(a)}{\leq} \sum_{i=1}^{j^*} \Eq{a_i(k)}
		\stackrel{(b)}{\leq} \sum_{i=1}^n \Eq{a_i(k)}
		\stackrel{(c)}{=}\lambda,
		\end{align*}
		where $(a)$ holds because $u_i(k)\leq s_i(k)$ with probability 1, by definition of unused service; $(b)$ holds because the number of arrivals to each queue is nonnegative; and $(c)$ holds by definition of the routing algorithm and $\lambda$.
	\end{or-letters}		
	This proves the theorem.
\end{pf}

\section{Heavy-traffic analysis}\label{sec:heavy-traffic}

In this section we perform heavy-traffic analysis of an heterogeneous load balancing system operating under \newline power-of-$d$ choices. Specifically, we prove that in the heavy-traffic limit, the load balancing system operating under power-of-$d$ choices behaves as a single server queue and we find the limiting joint distribution of the vector of queue lengths scaled by the heavy-traffic parameter.

Heavy traffic means that we load the system close to its maximum capacity. To take the limit we parametrize the system as follows. Fix a sequence of service rate vectors $\{\vs(k):k\in \bZ_+\}$ and take $\epsilon\in\left(0,\mu_\Sigma\right)$. The arrival process to the system parametrized by $\epsilon$ is an i.i.d. sequence $\{a^\peps(k):k\in \bZ_+\}$ that satisfies $\lambda^\peps\defn\E{a^\peps(1)}=\mu_\Sigma-\epsilon$. Then, the heavy-traffic limit is obtained by taking $\epsilon\dto 0$. 
We add a superscript $\peps$ to the queue length, arrival and unused service variables when we refer to the load balancing system parametrized by $\epsilon$.

In the next proposition we show SSC to a one-dimensional subspace. In other words, we show that, in the limit, the $n$-dimensional load balancing system operating under power-of-$d$ choices behaves as a single server queue. Before showing the result we introduce the following notation. For any vector $\vx\in\bR^n$, define
\begin{align}\label{eq:q-parallel-and-perp}
\vx_\parallel = \vone \left(\dfrac{\sum_{i=1}^n x_i}{n}\right)\;,\; \vx_\perp\defn \vx-\vx_\parallel.
\end{align}
Then, $\vx_\parallel$ is the projection of $\vx$ on the line $\{\vy\in\bR^n:y_i=y_j\,\forall i,j\in[n]\}$ and $\vx_\perp$ is the error of approximating $\vx$ by $\vx_\parallel$. Now we present the result.

\begin{proposition}\label{prop:SSC}
	Given a sequence $\{\vs(k):k\in \bZ_+\}$ of i.i.d. random vectors, and $\epsilon\in(0,\mu_\Sigma)$, consider a load balancing system operating under power-of-$d$ choices, \newline 
	parametrized by $\epsilon$ as described above. Suppose $d\geq 2$, and that 
	the number of arrivals and the potential service in each time slot are bounded. Let $\vmu\in\interior{\cM^\pd}$ and let $\vqbar^\peps$ be a steady-state vector such that $\{\vq^\peps(k):k\in \bZ_+\}$ converges in distribution to $\vqbar^\peps$ as $k\to\infty$.
	Let $\delta>0$ be such that for all $j\in\bZ_+$ satisfying $d\leq j\leq n-1$ we have
	\begin{align}\label{eq:beta-partial-sums-partial1}
	\dfrac{\sum_{i=1}^{j}\mu_{(i)}}{\mu_\Sigma}-\delta\geq \dfrac{\binom{j}{d}}{\binom{n}{d}}
	\end{align}
	If $\epsilon<\delta\mu_\Sigma$, then $\E{\|\vqbar_\perp^\peps \|^m}\leq M_m$ for each $m=1,2,\ldots$, where $M_m$ is a finite constant (independent of $\epsilon$).
\end{proposition}
Proposition \ref{prop:SSC} says that the error of approximating $\vqbar^\peps$ by $\vqbar_\parallel^\peps$ is negligible in heavy traffic because, as $\epsilon$ gets smaller, the arrival rate to the system increases and, therefore, the vector of queue lengths $\vqbar^\peps$ becomes larger. Then, the projection $\vqbar_\parallel^\peps$ also becomes larger. However, the error of approximating $\vqbar^\peps$ by $\vqbar_\parallel^\peps$, denoted as $\vqbar_\perp^\peps$, has bounded moments. Then, as $\epsilon$ goes to zero it becomes negligible.

Observe that the vector $\vqbar^\peps$ is well defined, because $\vmu\in\interior{\cM^\pd}\subset\cM^\pd$. Then, from Theorem \ref{thm:throughput-optimality}, for all $\epsilon>0$ the DTMC $\left\{\vq^\peps(k):k\in \bZ_+\right\}$ is positive recurrent.

In the proof of Proposition \ref{prop:SSC} we use a result first presented in \cite[Lemma 1]{atilla}, which is a corollary of a result proved in \cite{hajek_drift}. 
\begin{techreport}
	We restate this result in \ref{app:hajek-lemma} for completeness.
\end{techreport}

\begin{pf}[of Proposition \ref{prop:SSC}]
	For ease of exposition, we omit the dependence on $\epsilon$ on the variables. Define $V(\vq)\defn\left\|\vq\right\|^2$, $V_\parallel(\vq)\defn\|\vq_\parallel\|^2$, $W_\perp(\vq)\defn \|\vq_\perp\|$.	
	\begin{techreport}
		To prove the Proposition we use Lemma \ref{lemma:hajek} with $Z(\vq)=W_\perp(\vq)$.
	\end{techreport}
	\begin{or-letters}
		We use the Lyapunov function  $W_\perp(\vq)$.
	\end{or-letters}
	We start with a fact first used in \cite{atilla}. Observe that $\|\vq_\perp\|=\sqrt{\|\vq_\perp\|^2}$ by definition of square root, and $f(x)=\sqrt{x}$ is a concave function. Then, by definition of concavity and the Pythagoras theorem, 
	\begin{align}\label{eq:concavity}
	\Delta W_\perp(\vq)\leq \dfrac{1}{2\left\|\vq_\perp \right\|}\left(\Delta V(\vq)-\Delta V_\parallel(\vq) \right).
	\end{align}
	
	\begin{techreport}
		Then,  to prove condition (C1), it suffices to upper bound $\Eq{\Delta V(\vq)}$ and lower bound $\Eq{\Delta V_\parallel(\vq)}$. We start with $\Eq{\Delta V(\vq)}$. 
	\end{techreport}
	\begin{or-letters}
		We need to show two conditions. In the first condition we show that $\Eq{\Delta W_\perp(\vq(k))}$ is negative if $\vq$ lies outside a bounded set, and in the second condition we show that $\Eq{\Delta W_\perp(\vq(k))}$ is bounded. 
		
		To prove the first one, we find an upper bound to $\Eq{\Delta V(\vq)}$ and a lower bound to $\Eq{\Delta V_\parallel(\vq)}$. We start with $\Eq{\Delta V(\vq)}$. 
	\end{or-letters}
	From the proof of Theorem \ref{thm:throughput-optimality}, we know \eqref{eq:foster-lyapunov-partial1} is satisfied. 
	We analyze the last term differently here. Defining $\phi(i)$ as in the proof of Theorem \ref{thm:throughput-optimality}, we have
	\begin{techreport}
		\begin{align*}
		& \Eq{\langle\vq,\va(k)-\vs(k)\rangle} \\
		=& \lambda\sum_{i=1}^{n-d+1} q_{(i)}\dfrac{\binom{n-i}{d-1}}{\binom{n}{d}} - \sum_{i=1}^n q_{(i)}\mu_{\phi(i)} \\
		\stackrel{(a)}{=}& -\epsilon\left(\dfrac{\sum_{i=1}^n q_i}{n}\right)+ \sum_{i=1}^{n-d+1} q_{(i)}\dfrac{\lambda\binom{n-i}{d-1}}{\binom{n}{d}} + \sum_{i=1}^n q_{(i)}\left(\dfrac{\epsilon}{n}-\mu_{\phi(i)} \right) \\
		\stackrel{(b)}{=}& -\epsilon\left(\dfrac{\sum_{i=1}^n q_i}{n}\right) + \sum_{i=1}^n q_{(i)}\beta_i
		\end{align*}
		where $(a)$ holds by adding and subtracting $\frac{\epsilon}{n}\left(\sum_{i=1}^n q_i\right)$, and reorganizing terms; and $(b)$ holds defining for each $i\in[n]$
		\begin{align}\label{eq:def-beta}
		\beta_i\defn \begin{cases}
		\dfrac{\binom{n-i}{d-1}}{\binom{n}{d}}\lambda +\dfrac{\epsilon}{n} -\mu_{\phi(i)} & \text{, if }1\leq i\leq n-d+1 \\
		\dfrac{\epsilon}{n} -\mu_{\phi(i)} & \text{, if }n-d+1< i\leq n
		\end{cases} 
		\end{align}
		Observe $\beta_i=\alpha_i+\frac{\epsilon}{n}$ for each $i\in[n]$, where $\alpha_i$ is defined in \eqref{eq:def-alpha}.
	\end{techreport}
	\begin{or-letters}
		\begin{align*}
		& \Eq{\langle\vq,\va(k)-\vs(k)\rangle} 
		 \stackrel{(a)}{=} 
		-\epsilon\left(\dfrac{\sum_{i=1}^n q_i}{n}\right) + \sum_{i=1}^n q_{(i)}\beta_i,
		\end{align*}
		where $(a)$ holds reorganizing terms, and defining
		\begin{align}\label{eq:def-beta}
		\beta_i\defn \begin{cases}
		\dfrac{\binom{n-i}{d-1}}{\binom{n}{d}}\lambda +\dfrac{\epsilon}{n} -\mu_{\phi(i)} & \text{, if }1\leq i\leq n-d+1 \\
		\dfrac{\epsilon}{n} -\mu_{\phi(i)} & \text{, if }n-d+1< i\leq n.
		\end{cases} 
		\end{align}
	\end{or-letters}
	
	\begin{claim}\label{claim:betas}
		The parameters $\beta_i$ defined in \eqref{eq:def-beta} satisfy 
		\begin{enumerate}
			\item\label{beta:n} $\beta_n\leq -\mu_{(1)}+\frac{\epsilon}{n}$.
			
			\item\label{beta:total-sum} $\ds \sum_{i=1}^n \beta_i=0$.
			
			\item\label{beta:partial-sum} For any $j\in\bZ_+$ satisfying $2\leq j\leq n-1$ we have $\ds  \sum_{i=j}^n \beta_i \leq -\delta\mu_\Sigma + \epsilon$.
		\end{enumerate}
	\end{claim}
	
	We prove Claim \ref{claim:betas} in Section \ref{app:proof-betas}. {\blue Observe that if $d=1$, the second property is not satisfied.} Using Claim \ref{claim:betas} we obtain
	\begin{align}
	\sum_{i=1}^n q_{(i)}\beta_i =& q_{(1)}\sum_{i=1}^n \beta_i + \sum_{j=2}^n \left(\sum_{i=j}^n \beta_i\right)\left(q_{(j)}-q_{(j-1)}\right) \nonumber \\
	\leq& \left(-\delta\mu_\Sigma +\epsilon\right)\left(q_{(n)}-q_{(1)}\right). \label{eq:q-beta-qn-q1}
	\end{align}
	
	Observe that, by definition of $\vq_\perp$, we have
	\begin{align*}
	\left\|\vq_\perp \right\|^2 =& \sum_{i=1}^n \left(q_i-\dfrac{\sum_{j=1}^n q_j}{n} \right)
	\stackrel{(a)}{\leq} n\left(q_{(n)}-q_{(1)}\right),
	\end{align*}
	where $(a)$ holds because $q_i\leq q_{(n)}$ for all $i\in[n]$ and $\frac{1}{n}\sum_{j=1}^n q_j\geq q_{(1)}$ by definition of $q_{(1)}$ and $q_{(n)}$. Using this result in \eqref{eq:q-beta-qn-q1} we obtain that 
	\begin{align*}
	\sum_{i=1}^n q_{(i)}\beta_i\leq& \left(\dfrac{-\delta\mu_\Sigma+\epsilon}{\sqrt{n}}\right) \left\|\vq_\perp\right\|
	\leq \left(\dfrac{-\delta\mu_\Sigma+\epsilon_0}{\sqrt{n}}\right) \left\|\vq_\perp\right\|,
	\end{align*}
	for any $\epsilon_0\in\left(0, \delta \mu_\Sigma\right)$. Therefore,
	\begin{align}\label{eq:SSC-C1-V}
	\begin{aligned}
	& \Eq{\Delta V(\vq(k))} \\
	\leq& K_1 -2\epsilon\left(\dfrac{\sum_{i=1}^n q_i}{n}\right) +2 \left(\dfrac{-\delta\mu_\Sigma + \epsilon_0}{\sqrt{n}}\right) \left\|\vq_\perp\right\|.
	\end{aligned}
	\end{align}
	
	To lower bound $\Eq{\Delta V_\parallel(\vq)}$ we only use properties of the norm and the unused service. 
	We obtain
	\begin{align}\label{eq:SSC-C1-Vpar}
	\Eq{\Delta V_{\parallel}(\vq(k))}\geq -2\epsilon\left(\dfrac{\sum_{i=1}^n q_i}{n}\right)-K_3,
	\end{align}
	where $K_3\defn 2n\smax^2$, and $\smax$ is a finite constant such that $s_i(1)\leq \smax$ for all $i\in[n]$ with probability 1. Using \eqref{eq:SSC-C1-V} and \eqref{eq:SSC-C1-Vpar} in \eqref{eq:concavity} we obtain
	\begin{align*}
	\Eq{\Delta W_\perp(\vq(k))}\leq& \dfrac{K_1+K_3}{2\left\|\vq_\perp\right\|}+\left(\dfrac{-\delta\mu_\Sigma+\epsilon_0}{\sqrt{n}}\right),
	\end{align*}
	\begin{techreport}
		which satisfies condition (C1) for
		\begin{align*}
		\kappa = \left(\dfrac{K_1+K_3}{2}\right)\left(-\eta+\dfrac{\delta\mu_\Sigma-\epsilon_0}{\sqrt{n}} \right)^{-1}.
		\end{align*}
		
		Condition (C2) is trivially satisfied because potential service and arrivals in one time slot are bounded random variables.
	\end{techreport}
	\begin{or-letters}
		which proves the first condition of \cite[Lemma 1]{atilla}. The second condition is trivially satisfied because the arrival and service random variables are bounded.
	\end{or-letters}
\end{pf}

Using SSC, we can completely determine the behavior of the vector of queue lengths in heavy-traffic. In the next proposition we provide this result.

\begin{thm}\label{theo:distribution}
	Consider a set of load balancing systems operating under power-of-$d$ as described in Proposition \ref{prop:SSC}. 
	Let $\sigma_a^\peps$ be the standard deviation of $a^\peps(1)$ and assume $\sigma_a=\lim_{\epsilon\dto0}\sigma_a^\peps$. 
	Then, $\epsilon\vqbar^\peps\Longrightarrow \Upsilon\vone$ as $\epsilon\dto 0$, where $\Upsilon$ is an exponential random variable with mean $\dfrac{1}{2n}\left(\sigma_a^2+\vone^T \Sigma_s\vone\right)$, and $\Longrightarrow$ denotes convergence in distribution.
\end{thm}

\begin{rmk}
	In Proposition \ref{prop:SSC} and Theorem \ref{theo:distribution} we assume that the set $\cM^\pd$ has nonempty interior. This can be proved by observing that, for $d\geq 2$, a vector of homogeneous service rates $\vmu = c \vone$ (with $c>0$) satisfies all the inequalities in \eqref{eq:Md}, and none of them is tight. Then, such $\vmu=c\vone\in\interior{\cM^\pd}$. On the other hand, when $d=1$, the set $\cM^\pd$ only contains the  homogeneous service rate vectors, which has an empty interior. Then, our heavy-traffic results 
	are not applicable. This is consistent with the fact that random routing is not heavy-traffic optimal.
\end{rmk}

\begin{pf}[of Theorem \ref{theo:distribution}]
	We use the MGF method \cite{Hurtado_transform_method}, which is a two-step procedure to compute the joint distribution of the scaled vector of queue lengths in heavy traffic, in queueing systems that experience one-dimensional SSC. In fact, our theorem is a corollary of \cite[Theorem 2]{Hurtado_transform_method}. We only verify that three conditions are satisfied.
	
	We first verify that the routing algorithm is throughput optimal, which holds from Theorem \ref{thm:throughput-optimality} because we assume $\vmu\in\cM^\pd$. The second condition is SSC to a one-dimensional subspace, which is satisfied by Proposition \ref{prop:SSC}. 
	\begin{techreport}
		In fact, the authors in \cite{Hurtado_transform_method} require a weaker notion of State Space Collapse, which is trivially satisfied here, after proving Proposition \ref{prop:SSC}.
	\end{techreport}
	The last condition is existence of the MGF of $\epsilon\sum_{i=1}^n \qbar_i$, which we formalize in Claim \ref{claim:existence-MGF} and prove in Section \ref{app:existence-MGF}.
	
	\begin{claim}\label{claim:existence-MGF}
		For the load balancing system described in Theorem \ref{theo:distribution},  there exists $\Theta>0$ such that $\E{e^{\theta\epsilon\sum_{i=1}^n \qbar_i^\peps}}$ is finite for all $\theta\in[-\Theta,\Theta]$.
	\end{claim}
\end{pf}

{\blue
\section{Generalization to other routing policies}\label{sec:generalization}

In this section we generalize the sufficient conditions in Theorem \ref{thm:throughput-optimality} to a larger class of routing policies. Instead of using power-of-$d$ choices, suppose the router randomly selects an arbitrary subset of servers, and then the arrivals are routed to the server with the shortest queue among these. Let $\pi:2^{[n]}\to [0,1]$ be the probability mass function that governs the set of servers that are randomly selected in each time slot. We call $\cR^\pi$ the routing algorithm described above.

\begin{thm}\label{thm:generalization}
	Given $\pi:2^{[n]}\to [0,1]$, consider a load balancing system as described in Section \ref{sec:model}, operating under $\cR^\pi$. For each subset $\cS\subseteq[n]$, let $\pi(\cS)$ be the probability of sampling the servers in the set $\cS$. Let $\cP\left([n]\right)$ be the set of permutations of the elements of the set $[n]$, and for each $\tau\in \cP([n])$ define
	\begin{align*}
		\cM_\tau \defn & \Bigg\{\vmu\in\bR^n_+: \dfrac{\ds\sum_{i=1}^j \mu_{(i)}}{\mu_\Sigma}\leq \sum_{i=1}^j \sum_{\cS\in \cS^\tau_i } \pi(\cS)\;\forall j\in[n-1] \Bigg\},
	\end{align*}
	where $\begin{array}[t]{ll}
			\cS^\tau_i \defn \left\{\cS\subseteq [n]:\right.& \tau(n-i+1)\in\cS,\\
		& \left.\tau(\ell)\notin \cS\quad \forall \ell<n-i+1 \right\}.
	\end{array}$

	\noindent $\cR^\pi$ is throughput optimal if $\vmu\in\cM_\tau$ for all $\tau\in \cP\left([n]\right)$.
\end{thm}

The proof is similar to the proof of Theorem \ref{thm:throughput-optimality}, and we present a sketch in Section \ref{app:generalization} for completeness.
}

\section{Details of the proofs in Sections \ref{sec:throughput-optimality}, \ref{sec:heavy-traffic} and \ref{sec:generalization}}\label{sec:proofs}

\subsection{Proof of Claim \ref{claim:alphas}}\label{app:proof-alphas}

\begin{techreport}
	Recall the definition of $\alpha_i$'s. For each $i\in[n]$ we have
	\begin{align*}
		\alpha_i\defn \begin{cases}
			\dfrac{\lambda\binom{n-i}{d-1}}{\binom{n}{d}} - \mu_{\phi(i)} & \text{, if }1\leq i\leq n-d+1 \\
			-\mu_{\phi(i)} & \text{, if }n-d+1< i\leq n
		\end{cases}
	\end{align*}
	
	Now we prove the claim.
\end{techreport}

\begin{pf}[of Claim \ref{claim:alphas}]
	
	We prove each of the three properties. We obtain:
	\begin{enumerate}
		\item If $i=n$ we have $\alpha_n = -\mu_{\phi(n)}\leq -\mu_{(1)}$, because $\mu_{(1)}=\min_{i\in[n]}\mu_i$.
		
		\item The total sum of $\alpha_i$'s satisfies
		\begin{align*}
			\sum_{i=1}^n\alpha_i =& \dfrac{\lambda}{\binom{n}{d}} \sum_{i=1}^{n-d+1} \binom{n-i}{d-1} - \mu_\Sigma
			\stackrel{(a)}{=} \lambda-\mu_\Sigma =-\epsilon,
		\end{align*}
		where $(a)$ holds because $\sum_{i=1}^{n-d+1} \binom{n-i}{d-1}=\binom{n}{d}$.
		
		\item If $2\leq j\leq n-d+1$  we have that the tail sums are
		\begin{techreport}
			\begin{align*}
				\sum_{i=j}^n \alpha_i=& \dfrac{\lambda}{\binom{n}{d}}\sum_{i=j}^{n-d+1} \binom{n-i}{d-1} - \sum_{i=j}^n \mu_{\phi(i)} \\
				\stackrel{(a)}{=}& \lambda \dfrac{\binom{n+1-j}{d}}{\binom{n}{d}} - \sum_{i=j}^n\mu_{\phi(i)} \\
				\stackrel{(b)}{=}& \dfrac{\binom{n+1-j}{d}}{\binom{n}{d}}(\mu_\Sigma-\epsilon) - \sum_{i=j}^n\mu_{\phi(i)} \\
				\stackrel{(c)}{\leq}& \sum_{i=1}^{n+1-j} \mu_{(i)}- \dfrac{\binom{n+1-j}{d}}{\binom{n}{d}}\epsilon - \sum_{i=j}^n\mu_{\phi(i)} \\
				\stackrel{(d)}{\leq}& -\dfrac{\epsilon}{\binom{n}{d}},
			\end{align*}
			where $(a)$ holds because $\sum_{i=j}^{n-d+1} \binom{n-i}{d-1}= \binom{n+1-j}{d}$; $(b)$ holds by definition of $\epsilon$; $(c)$ holds because $\vmu\in\cM^\pd$; and $(d)$ holds because $\binom{n+1-j}{d}\geq 1$, and because  $\sum_{i=1}^{n+1-j}\mu_{(i)}- \sum_{i=j}^n \mu_{\phi(i)} \leq 0$, since $\sum_{i=1}^{n+1-j}\mu_{(i)}$ is the sum of the $n+j-1$ smallest elements of $\vmu$, and $\sum_{i=j}^n \mu_{\phi(i)}$ is the sum of $n+j-1$ of the elements of $\vmu$ which are not necessarily the smallest.
		\end{techreport}
		\begin{or-letters}
			\begin{align*}
				\sum_{i=j}^n \alpha_i
				\stackrel{(a)}{=}& \lambda \dfrac{\binom{n+1-j}{d}}{\binom{n}{d}} - \sum_{i=j}^n\mu_{\phi(i)} \\
				\stackrel{(b)}{\leq}& \sum_{i=1}^{n+1-j} \mu_i- \dfrac{\binom{n+1-j}{d}}{\binom{n}{d}}\epsilon - \sum_{i=j}^n\mu_{\phi(i)}
				\stackrel{(c)}{\leq} -\dfrac{\epsilon}{\binom{n}{d}},
			\end{align*}
			where $(a)$ holds because $\sum_{i=j}^{n-d+1} \binom{n-i}{d-1}= \binom{n+1-j}{d}$; $(b)$ holds by definition of $\epsilon$ and because $\vmu\in\cM^\pd$; and $(c)$ holds because $\binom{n+1-j}{d}\geq 1$, and because $\sum_{i=1}^{n+1-j}\mu_i$ is the sum of the $n+j-1$ smallest elements of $\vmu$.
		\end{or-letters}		
		If $n-d+1<j\leq n-1$ we have
		\begin{align*}
			\sum_{i=j}^n \alpha_i=& -\sum_{i=j}^n \mu_{\phi(i)} \leq -\mu_{(1)}.
		\end{align*}
		\begin{techreport}
			Then, for all $2\leq j\leq n-1$ we have
			\begin{align*}
				\sum_{i=j}^n \alpha_i \leq -K_2\defn-\min\left\{\dfrac{\epsilon}{\binom{n}{d}},\, \mu_{(1)} \right\}.
			\end{align*}
		\end{techreport}
	\end{enumerate}

\end{pf}

\subsection{Proof of Claim \ref{claim:betas}}\label{app:proof-betas}

\begin{techreport}
	Recall the definition of $\beta_i's$. For each $i\in[n]$ we have
	\begin{align*}
		\beta_i\defn \begin{cases}
			\dfrac{\binom{n-i}{d-1}}{\binom{n}{d}}\lambda +\dfrac{\epsilon}{n} -\mu_{\phi(i)} & \text{, if }1\leq i\leq n-d+1 \\
			\dfrac{\epsilon}{n} -\mu_{\phi(i)} & \text{, if }n-d+1< i\leq n
		\end{cases} 
	\end{align*}
\end{techreport}

\begin{pf}[of Claim \ref{claim:betas}]
	\begin{techreport}
		We prove each of the three properties. We have:
		\begin{enumerate}
			\item If $i=n$ we have
			\begin{align*}
				\beta_n=& \alpha_n+\frac{\epsilon}{n}\leq -\mu_{(1)}+\frac{\epsilon}{n},
			\end{align*}
			where we used property \ref{alpha:prop-n} from Claim \ref{claim:alphas}.
			
			\item The total sum of $\beta_i$'s satisfies
			\begin{align*}
				\sum_{i=1}^n \beta_i =& \sum_{i=1}^n \alpha_i  + \epsilon = 0,
			\end{align*}
			where we used property \ref{alpha:total-sum} from Claim \ref{claim:alphas}.
			
			\item To prove this property we divide in 2 cases. If $j\leq n-d+1$ we have
			\begin{align*}
				\sum_{i=j}^n \beta_i =& \sum_{i=j}^{n-d+1} \dfrac{\binom{n-i}{d-1}}{\binom{n}{d}}\lambda + \sum_{i=j}^n \left(\dfrac{\epsilon}{n}-\mu_{\phi(i)}  \right) \\
				=& \dfrac{\binom{n+1-j}{d}}{\binom{n}{d}}\left(\mu_\Sigma -\epsilon \right) + \dfrac{n-j+1}{n}\epsilon - \sum_{i=j}^{n}\mu_{\phi(i)} \\
				\stackrel{(a)}{\leq}& \dfrac{\binom{n+1-j}{d}}{\binom{n}{d}}\mu_{\Sigma} + \epsilon - \sum_{i=1}^{n-j+1}\mu_{(i)} \\
				\stackrel{(b)}{\leq}& \epsilon -\delta\mu_{\sigma}
			\end{align*}
			where $(a)$ holds because $\epsilon>0$, $\frac{n-j+1}{n}\leq 1$ and because $\sum_{i=1}^{n-j+1}\mu_{(i)}$ is the sum of the smallest $(n-j+1)$ elements of $\vmu$; and $(b)$ holds by \eqref{eq:beta-partial-sums-partial1} and reorganizing terms.
			
			If $j>n-d+1$ we have
			\begin{align*}
				\sum_{i=j}^n \beta_i=& \sum_{i=j}^n \left(\dfrac{\epsilon}{n}-\mu_{\phi(i)}  \right) \\
				\stackrel{(a)}{\leq}& \dfrac{n-j+1}{n}\epsilon - \sum_{i=1}^{n-j+1}\mu_{(i)} \\
				\stackrel{(b)}{\leq}& \epsilon - \mu_\Sigma\left(\dfrac{\binom{n-j+1}{d}}{\binom{n}{d}} + \delta \right) \\
				\stackrel{(c)}{\leq}& \epsilon -\delta\mu_\Sigma
			\end{align*}
			where $(a)$ holds because $\sum_{i=1}^{n-j+1}\mu_{(i)}$ is the sum of the smallest $(n-j+1)$ elements of $\vmu$; $(b)$ holds because $\frac{n-j+1}{n}\leq 1$ and by \eqref{eq:beta-partial-sums-partial1}; and $(c)$ because $\frac{\binom{n-j+1}{d}}{\binom{n}{d}}\geq 0$.
		\end{enumerate}
	\end{techreport}
	\begin{or-letters}
		Properties 1 and 2 follow immediately from the fact that $\beta_i=\alpha_i+\frac{\epsilon}{n}$. To prove the third property we divide in cases. If $j\leq n-d+1$ we have
		\begin{align*}
			\sum_{i=j}^n \beta_i &= \sum_{i=j}^{n-d+1} \dfrac{\binom{n-i}{d-1}}{\binom{n}{d}}\lambda + \sum_{i=j}^n \left(\dfrac{\epsilon}{n}-\mu_{\phi(i)}  \right) \\
			&\leq  \dfrac{\binom{n+1-j}{d}}{\binom{n}{d}}\mu_{\Sigma} + \epsilon - \sum_{i=1}^{n-j+1}\mu_{i} 
			\stackrel{(a)}{\leq} \epsilon -\delta\mu_{\Sigma},
		\end{align*}
		where $(a)$ 
		holds by \eqref{eq:beta-partial-sums-partial1} and reorganizing terms.
		
		If $j>n-d+1$ we have
		\begin{align*}
			\sum_{i=j}^n \beta_i &= \sum_{i=j}^n \left(\dfrac{\epsilon}{n}-\mu_{\phi(i)}  \right) \\
			&\leq  \dfrac{n-j+1}{n}\epsilon - \sum_{i=1}^{n-j+1}\mu_{i}
			\stackrel{(a)}{\leq} \epsilon -\delta\mu_\Sigma
		\end{align*}
		where $(a)$ holds 
		by \eqref{eq:beta-partial-sums-partial1} and because $\frac{n-j+1}{n}\leq 1$ .
	\end{or-letters}
	
\end{pf}

\begin{techreport}
	\section{Preliminary results for the proof of Theorem \ref{thm:throughput-optimality}}\label{app:F-L.theorem}
	
	We first present Foster-Lyapunov theorem as stated in \cite[Theorem 3.3.7]{srikantleibook}.
	
	\begin{thm}\label{lemma:foster-lyapunov}
		Let $\{X(k):k\in \bZ_+\}$ be an irreducible Markov chain with state space $\cS$. Suppose that there exists a function $V:\cS\to\bR_+$ and a finite set $\cB\subseteq\cS$ satisfying the conditions
		\begin{itemize}
			\item[(C1)] $\E{V(X(k+1))-V(X(k)) | X_k=x }\leq -\xi$ if $x\in\cS\setminus\cB$ for some $\xi>0$
			\item[(C2)] $\E{V(X(k+1))-V(X(k)) | X_k=x }\leq \kappa$ if $x\in\cB$ for some $\kappa<\infty$
		\end{itemize}
		Then, the Markov chain $\left\{X(k):k\in \bZ_+ \right\}$ is positive recurrent.
	\end{thm}
	
	Now we present a certificate that a Markov chain is not positive recurrent \cite[Theorem 3.3.10]{srikantleibook}.
	
	\begin{lem}\label{lemma:foster-lyapunov-not-pos-rec}
		An irreducible Markov chain $\{X(k):k\in \bZ_+ \}$ with state space $\cS$ is not positive recurrent (i.e., it is either transient or null recurrent) if there exists a function $V:\cS\to\bR_+$ and a finite set $\cB\subseteq \cS$ satisfying the following conditions
		\begin{itemize}
			\item[(C1)] $\E{V(X(k+1))-V(X(k))|X(k)=x}\geq 0$ for all $x\in\cS\setminus \cB$
			\item[(C2)] There exists some $x\in\cS\setminus \cB$ such that $V(x)>V(y)$ for all $y\in\cB$
			\item[(C3)] $\E{|V(X(k+1))-V(X(k))| \,|X(k)=x}\leq \kappa$ for all $x\in\cS$ and some $\kappa<\infty$
		\end{itemize}
	\end{lem}

	\subsection{Preliminary result for the proof of Proposition \ref{prop:SSC}}\label{app:hajek-lemma}
	
	\begin{lem}\label{lemma:hajek}
		For an irreducible and aperiodic Markov Chain $\{X(k):\,k\in \bZ_+\}$ over a countable state space $\cS$, suppose $Z:\cS\to \bR_+$ is a nonnegative valued Lyapunov function. The drift of $Z$ at $x$ is
		\begin{align*}
			\Delta Z(x)\defn\big[Z\big(X(k+1)\big)-Z\big(X(k)\big) \big]\ind{X(k)=x}
		\end{align*}
		Thus, $\Delta Z(x)$ is a random variable that measures the amount of change in the value of $Z$ in one step, starting from state $x$. This drift is assumed to satisfy the following conditions:
		\begin{itemize}
			\item[(C1)] There exists $\eta>0$ and $\kappa<\infty$ such that
			\begin{align*}
				\E{\left.\Delta Z(x)\,\right|\,X(k)=x}\leq -\eta\quad\text{for all $x\in\cS$ with $Z(x)\geq \kappa$}
			\end{align*}
			
			\item[(C2)] There exists $D<\infty$ such that
			\begin{align*}
				|\Delta Z(x)|\leq D\quad\text{ with probability 1 for all $x\in\cS$}
			\end{align*}
		\end{itemize}
		If we further assume that the Markov chain $\{X(k):\,k\in \bZ_+\}$ is positive recurrent, then $Z(X(k))$ converges in distribution to a random variable $\overline{Z}$ for which
		\begin{align*}
			\E{e^{\theta^* \overline{Z}}}\leq C^*
		\end{align*}
	\end{lem}
	
\end{techreport}

\subsection{Existence of MGF}\label{app:existence-MGF}

\begin{pf}[of Claim \ref{claim:existence-MGF}]
	The proof is similar to the proof of existence of MGF under JSQ routing, which was done in \cite{Hurtado_transform_method}. We write a sketch of the proof here for completeness. 
	First observe that if $\theta\leq 0$, the proof holds trivially. 
	Now, assume $\theta>0$. Observe that $f(x)=e^x$ is a convex function. Then, by Jensen's inequality, 
	we have
	\begin{align*}
		e^{\frac{\theta}{n}\epsilon\sum_{i=1}^n q_i}\leq \dfrac{1}{n}\sum_{i=1}^n e^{\theta\epsilon q_i}.
	\end{align*}
	Then, it suffices to show that $\E{\sum_{i=1}^n e^{\theta\epsilon q_i}}<\infty$ for $\theta\in[-\Theta,\Theta]$, for all $i\in[n]$. 
	\begin{techreport}
		We use Theorem \ref{lemma:foster-lyapunov} with function $V_{MGF}(\vq)=\sum_{i=1}^n e^{\theta\epsilon q_i}$.
	\end{techreport}
	\begin{or-letters}
		We use Foster-Lyapunov theorem \cite[Theorem 3.3.7]{srikantleibook} with function $V_{MGF}(\vq)=\sum_{i=1}^n e^{\theta\epsilon q_i}$.
	\end{or-letters}
	For each $i\in[n]$ we have
	\begin{align*}
		\left(e^{\theta\epsilon q_i(k+1)}-1 \right)\left(e^{-\theta\epsilon u_i(k)}-1\right)=0,
	\end{align*}
	which holds by \eqref{eq:qu}. 
	\begin{techreport}
		Then, reorganizing terms we have
		\begin{align*}
			e^{\theta\epsilon q_i(k+1)}= 1-e^{-\theta \epsilon u_i(k)}+ e^{\theta\epsilon \left(q_i(k)+a_i(k)-s_i(k)\right)}.
		\end{align*}
		Then, we obtain
	\end{techreport}
	\begin{or-letters}
		Then, reorganizing terms and summing over $i\in[n]$ we have
	\end{or-letters}
	\begin{align}
		& \Eq{\Delta V_{MGF}(\vq(k))} \nonumber \\
		& \begin{aligned}\label{eq:mgf-drift}
			&= \sum_{i=1}^n \left(1-\E{e^{-\theta\epsilon u_i(k)}}\right) \\ 
			& + \sum_{i=1}^n e^{\theta\epsilon q_{(i)}}\left(\Eq{e^{\theta\epsilon \left(a_{\phi(i)}(k)-s_{\phi(i)}(k) \right)}} -1\right),
		\end{aligned}
	\end{align}
	where $\phi(i)$ is defined as in the proof of Theorem \ref{thm:throughput-optimality}. 
	Since $\vu(k)\geq \vzero$ and $\theta>0$, the first term is upper bounded by $n$. Now, for a bounded random variable $X$, define $M_X(\theta)\defn \E{e^{\theta\epsilon X}}$. Then, for each $i\in[n]$ we have
	\begin{align*}
		& \Eq{e^{\theta\epsilon\left(a_{\phi(i)}(k)-s_{\phi(i)}(k) \right)}} -1 \\
		=& M_{a_{\phi(i)}-s_{\phi(i)}}(\theta)-1 
		\stackrel{(a)}{=} \theta M'_{a_{\phi(i)}-s_{\phi(i)}}(\xi_i),
	\end{align*}
	where $(a)$ holds by Taylor expansion, for a number $\xi_i$ between $0$ and $\theta$. Observe that the MGF is continuously differentiable at $\theta=0$ \cite[p.78]{mood} and
	\begin{align*}
		M_{a_{\phi(i)}-s_{\phi(i)}}'(0)=& \Eq{a_{\phi(i)}(k)-s_{\phi(i)}(k)}
		= \alpha_i,
	\end{align*}
	where $\alpha_i$ was defined in \eqref{eq:def-alpha}. For each $i\in[n]$, let $\tilde{\Theta}_i>0$ be such that for all $\theta$ between 0 and $\tilde{\Theta}_i$ we have
	\begin{align*}
		M_{a_{\phi(i)}-s_{\phi(i)}}'(\xi_i)\leq \frac{1}{2}\alpha_i.
	\end{align*}
	Let $\tilde{\Theta}=\min_{i\in[n]}\tilde{\Theta}_i$. Then, for all $\theta$ satisfying $\theta\epsilon<\tilde{\Theta}$
	\begin{align*}
		& \sum_{i=1}^n e^{\theta\epsilon q_{(i)}} \left(\Eq{e^{\theta\epsilon\left(a_{\phi(i)}(k)-s_{\phi(i)}(k) \right)}} -1\right) 
		\leq \sum_{i=1}^n e^{\theta\epsilon q_{(i)}}\alpha_i. 
	\end{align*}
	The rest of the proof is equivalent to the last steps of the proof of throughput optimality, so we omit it for brevity. The proof concludes by letting $\Theta=n\tilde{\Theta}$.
\end{pf}

{\blue
	\subsection{Proof of Theorem \ref{thm:generalization}}\label{app:generalization}
	
	\begin{pf}[of Theorem \ref{thm:generalization}]
		The proof is very similar to Theorem \ref{thm:throughput-optimality}. In fact, the only difference is the computation of $\Eq{\langle\vq,\va(k)\rangle}$. Since the sampling scheme in power-of-$d$ choices is symmetric, in Theorem \ref{thm:throughput-optimality} we obtain the simple expression presented in \eqref{eq:Ea}. In this case, we obtain
		\begin{align*}
			\Eq{\langle \vq,\va(k)\rangle} &= \sum_{i=1}^{n}q_{(i)}\lambda \Bigg(\sum_{\substack{ \cS\subseteq [n]: \\ \phi(i)\in\argmin_{\ell\in\cS} q_\ell }} \pi(\cS) \Bigg).
		\end{align*}
		We omit the rest of the proof for brevity.
	\end{pf}
}

\section{Acknowledgments}

This work was partially supported by the National Science Foundation [NSF-CCF: 1850439].
Daniela Hurtado-Lange has partial funding from ANID/DOCTORADO BECAS CHILE/2018 [72190413]

\bibliography{../../biblio-ok}


\end{document}